\documentclass[12pt,twoside]{article}

\usepackage{graphicx}
\usepackage{amssymb,amsmath}

\newcommand{\T}{{\mathcal T}}
\newcommand{\Tb}{{\mathbb T}}

\newcommand{\N}{{\mathbb N}}
\newcommand{\E}{{\mathbb E}}

\renewcommand{\P}{{\mathbb P}}

\renewcommand{\P}{{\mathbb P}}

\newcommand{\0}{{\mathbf 0}}

\newtheorem{stelling}{Theorem}

\numberwithin{equation}{section}

\begin{document}




\author{\textbf{Eric Cator, Henk Don}\\
\\
\footnotesize Institute for Mathematics, Astrophysics, and Particle Physics\\
\footnotesize Faculty of Science, Radboud University Nijmegen; e-mail:e.cator@math.ru.nl, henkdon@gmail.com}

\title{Constructing and searching conditioned Galton-Watson trees}
\maketitle

\begin{abstract}
We investigate conditioning Galton-Watson trees on general recursive-type events, such as the event that the tree survives until a specific level. It turns out that the conditioned tree is again a type of Galton-Watson tree, with different types of offspring and a level-dependent offspring distribution, which will all be given explicitly. As an interesting application of these results, we will calculate the expected cost of searching a tree until reaching a given level.
\end{abstract}
%

\section{Introduction}

The asymptotic shape of conditioned Galton-Watson trees has been
widely studied. For example, one could condition on the number of
nodes of the tree being $n$, and letting $n\rightarrow\infty$. In
a lot of cases, the limiting tree is quite well understood, see
for example the survey paper by Janson \cite{Janson}. Some work on
finite conditioned trees has been done by Geiger and Kersting \cite{Geiger},
who studies the shape of a tree conditioned on having height
exactly equal to $n$.

In this paper, we will investigate conditioning on events of a
recursive nature (as explained in section
\ref{section:generalization}), one example being conditioning on
survival to a given level. In sections \ref{section:condsurv} and
\ref{section:proof} this example will be worked out. The main idea
is that we consider different types of nodes, and the type of a node is determined by the type of her children. The offspring distribution of a node depends on its type and
on the level of the tree where this node is living. For trees conditioned to survive until a given level, we give a direct construction in Section
\ref{section:condsurv}. In Section \ref{section:proof} we show
that this construction indeed leads to the desired conditioned
tree.

Section \ref{section:treesearch} gives an application: determining
the cost of searching a tree to a given level. The model takes
into account costs for having a lot of children, but also for
walking into dead ends. The results from sections
\ref{section:condsurv} and \ref{section:proof} lead to recursions
that enable us to calculate the costs. For the case of Poisson
offspring these recursions are analyzed, including an
investigation of the asymptotically optimal mean offspring.

Finally in Section \ref{section:generalization} we give a
general setup that can be used to condition trees on suitable
events. Recursive relations that characterize such events will
be given, together with some illustrative examples. In particular, we will show how conditioning a tree on having height exactly $n$ fits into this framework. Our results therefore contain the case that has been studied by Geiger and Kersting.   

\subsection{Notation and preliminaries}

We will consider rooted Galton-Watson trees with arbitrary offspring
distribution, that can depend on the current height (or generation). We use the notation $\0$ to denote the root of a
tree. Define the set of trees of height $0$ as $\T_0=\{\0\}$, so
only the tree consisting of the root. Then we define inductively
for $k\geq 1$ the set of trees of height $k$:
\[ \T_k = \{ \0\}\sqcup \bigsqcup_{n=1}^\infty \{ [T_1,\ldots, T_n]\mid T_i\in \T_{k-1}\}.\]
For a tree $T = [T_1,\ldots, T_n] \in \T_{k-1}^n$, the trees
$T_1,\ldots, T_n$ will be called the children of $T$. We now
define the Galton-Watson probability measure on $\T_k$, with
offspring distribution $\mu_l$ at height $l\geq 0$, for arbitrary probability measures
$\mu_l$ on $\N(=\{0,1,\ldots\})$. Define independent random variables $W_l\sim
\mu_l$. Define $P_0$ as the trivial probability measure on $\T_0$: 
$P_0(\0)=1$. Now define inductively the probability measure $P_{lk}$ for $0\leq l\leq k$ and $k\geq 1$ as the following
probability measure on $\T_{k-l}$: if $l=k$, then $P_{kk}=P_0$, otherwise
\[ P_{lk}(\0) = \P(W_{l}=0)\]
and for all $T_1,\ldots, T_n\in \T_{k-l-1}$
\[P_{lk}([T_1,\ldots,T_n]) = \P(W_{l}=n)\prod_{i=1}^n P_{l+1,k}(T_i).\]
The intuition is that the second index determines the size of the final tree we are considering, whereas the first index determines at which level we are building up the tree (so $P_{lk}$ generates trees at level $l$ of size $k-l$). We are interested in $P_{0k}$, which is the Galton-Watson probability measure on $\T_k$ (trees cut off at height $k$), with offspring distribution $\mu_l$ on height $l$.

\section{Conditioning on survival to level $k$}\label{section:condsurv}

In this section we will define an alternative measure
$\tilde{P}_{lk}$ on $\T_{k-l}$, that allows us to efficiently construct a
Galton-Watson tree conditioned on reaching level $k$ (or on
extinction before level $k$). In the next section we prove that
indeed $P_{lk}=\tilde{P}_{lk}$. We fix $k\geq 1$ and for $0\leq l\leq k$
we define $A_l\subset \T_l$ to be the set of all trees $T\in\T_l$
that reach level $l$ ($A_0=\T_0$). Since $T\in A_l$ if and only if at least one
of its children reaches level $l-1$, the events $A_l$ satisfy the
following recursive relation
\begin{equation}\label{eq:recursionA_i}
A_0 = \T_0 \quad\textrm{and}\quad A_l = \left\{[T_1,\ldots,
T_n]\in \T_{l-1}^n : \#\left\{j:T_j\in A_{l-1}\right\}\geq
1\right\}.
\end{equation}
Let $p_{lk}$ be the probability that a tree starting at level $l$ reaches level $k$, so for $0\leq l\leq k$
\begin{equation}\label{eq:p_lk}
p_{lk} = P_{lk}(A_{k-l}).
\end{equation}
We will construct $\tilde P_{lk}$ as a Galton-Watson tree with two
types of children. A child born on level $l\leq k$ is of type 1 if
it is an element of $A_{k-l}$ (that is, if it has at least one
descendant at level $k$ in the original tree). All other children
are of type 2. Consequently, our tree will have different
offspring distribution for each level.

Introduce the independent random vectors $(W_{lk},X_{lk})\in\mathbb{N}^2, 0\leq
l\leq k-1$, where
$$
W_{lk}\sim \mu_{l} \quad{\rm and}\quad X_{lk}\mid W_{lk}\sim {\rm
Bin}(W_{lk},p_{l+1,k}).
$$
Here $W_{lk}$ represents the total number of children of an $l$th
level node and $X_{lk}$ the number of type 1 children. We start by
defining probability measures $\tilde Q_{lk}$ on $A_{k-l}\subset \T_{k-l}$,
for $0\leq l\leq k$, and probability measures $\tilde R_{lk}$ on
$A_{k-l}^c\subset \T_{k-l}$, for $1\leq l\leq k$. We take
$$
\tilde{Q}_{kk}(\0) = 1,
$$
and for $0\leq l\leq k-1$ we let
\begin{equation}\label{eq:root}
\tilde{Q}_{lk}(\0) = 0 \quad\textrm{and}\quad \tilde{R}_{lk}(\0) =
\P(W_{lk}=0 \mid X_{lk}=0).
\end{equation}
For all other trees the measures $\tilde Q_{lk}$ and $\tilde R_{lk}$
will be defined inductively. Let $1\leq l\leq k$ and choose
$T=[T_1,\ldots, T_n]\in \T_{k-l-1}^n$, for some $n\geq 1$. Define
$m=\#\left\{j:T_j\in A_{k-l-1}\right\}$. Then, if $m\geq 1$ (this is
equivalent to $T\in A_{k-l}$), define
\begin{equation}\label{eq:inAk}
\tilde{Q}_{lk}(T) = \frac{\P(W_{lk}=n, X_{lk}=m\mid X_{lk}\geq 1)}{{n \choose
m}} \prod_{i:T_i\in A_{k-l-1}} \tilde{Q}_{l+1,k}(T_i)\prod_{j:T_j\in
A_{k-l-1}^c} \tilde{R}_{l+1,k}(T_j),
\end{equation}
(defining products over empty sets to be 1) and
\[\tilde{R}_{lk}(T)=0.\]
If, on the other hand, $m=0$ (and therefore $T\not\in A_{k-l}$), we
define
\[ \tilde{Q}_{lk}(T) =0\]
and
\begin{equation}\label{eq:notinAk}
\tilde{R}_{lk}(T) =  \P(W_{lk}=n \mid X_{lk}=0) \prod_{j=1}^n
\tilde{R}_{l+1,k}(T_j).
\end{equation}
Note that if $l=k-1$, $T = [\0,\ldots,\0]$, so $m=n\geq 1$ and the
second product in (\ref{eq:inAk}) is empty. So $\tilde{R}_{kk}$ does
not play a role in the definition of $\tilde{Q}_{lk}$ and
$\tilde{R}_{lk}$ with $0\leq l\leq k-1$. Finally, we define for $0\leq
l\leq k$ (note that $p_{kk}=1$)
\begin{equation}\label{eq:Ptilde}
\tilde{P}_{lk}(T) = p_{lk} \tilde{Q}_{lk}(T) + (1-p_{lk}) \tilde{R}_{lk}(T)
\end{equation}
as a measure on $\T_{k-l}$.\\

We can describe the random tree $\Tb\sim \tilde{P}_{lk}$ as follows:
we first toss a coin to determine whether $\Tb$ is in $A_{k-l}$
(probability $p_{lk}$) or not (probability $1-p_{lk}$). If $\Tb\in A_{k-l}$,
then we choose it according to $\tilde{Q}_{lk}$, otherwise we choose
it according to $\tilde{R}_{lk}$. When $\Tb\sim \tilde{Q}_{lk}$, we
choose $(\tilde{W}_{lk},\tilde{X}_{lk})$, where $\tilde{W}_{lk}$ is the
number of children and $\tilde{X}_{lk}$ the number of children that
will lie in $A_{k-l-1}$, according to $(\tilde{W}_{lk},\tilde{X}_{lk})\sim
(W_{lk},X_{lk})\mid X_{lk}\geq 1$. The $\tilde{X}_{lk}$ children that lie in
$A_{k-l-1}$ are distributed over the $\tilde{W}_{lk}$ positions
uniformly at random. Then for each child in $A_{k-l-1}$ we draw a
tree according to $\tilde{Q}_{l+1,k}$, and for each child not in
$A_{k-l-1}$, we draw a tree according to $\tilde{R}_{l+1,k}$. If, on
the other hand, $\Tb \sim \tilde{R}_{lk}$, then we choose
$(\tilde{W}_{lk},\tilde{X}_{lk})\sim(W_{lk},X_{lk})\mid X_{lk}=0$. Since
$\tilde{X}_{lk}=0$, all the $\tilde{W}_{lk}$ children are in $A_{k-l-1}^c$
and for each of them, we draw a tree according to
$\tilde{R}_{l+1,k}$.

In this way we have described the random tree as a Galton-Watson
tree with two types of children and different offspring
distribution at each level. Note that conditioning $\tilde{P}_{lk}$
on $A_{k-l}$ is trivial: we simply have to draw $\Tb$ according to
$\tilde{Q}_{lk}$. A similar conclusion holds for conditioning on
$A_{k-l}^c$, in which case $\Tb\sim \tilde{R}_{lk}$.

\section{The two random trees are equally distributed}\label{section:proof}

Fix $k\geq 1$ and define all measures as before. We will check for
all $0\leq l\leq k$ and every $T\in \T_{k-l}$ that the probability of
that tree under $P_{lk}$ and $\tilde{P}_{lk}$ is the same:

\begin{stelling} For all $0\leq l\leq k$ and $T\in \T_{k-l}$,
$$
P_{lk}(T)=\tilde{P}_{lk}(T).
$$
\end{stelling}

\textbf{Proof.} Clearly $\tilde{P}_{kk}=P_{kk}=P_0$. We continue using
induction on $l$, where we fix $k$: let $0\leq l\leq k-1$, suppose $P_{l+1,k}=\tilde{P}_{l+1,k}$
and take $T\in \T_{k-l}$. If $T=\0$, by (\ref{eq:root}) and
(\ref{eq:Ptilde}),
$$
\tilde{P}_{lk}(\0) = (1-p_{lk})\P(W_{lk}=0 \mid X_{lk}=0) =
(1-p_{lk})\frac{\P(W_{lk}=0)}{\P(X_{lk}=0)},
$$
where the second equality follows from the fact that $W_{lk}=0$
implies $X_{lk}=0$. Furthermore
\begin{eqnarray*}
\P(X_{lk}=0) &=& \sum_{n=0}^\infty \P(W_{lk}=n) \P(X_{lk}=0\mid W_{lk}=n)\\
 &=& \sum_{n=0}^\infty \P(W_{lk}=n) \P({\rm Bin}(n,p_{l+1,k})=0) = \sum_{n=0}^\infty \P(W_{lk}=n) (1-p_{l+1,k})^n.
\end{eqnarray*}
This last expression is precisely the probability (under $P_{lk}$)
that all children of a tree starting at level $l$ are in $A_{k-l-1}^c$. This is by (\ref{eq:recursionA_i}) exactly the case if the tree itself is in $A_{k-l}^c$, implying that  
$$
\P(X_{lk}=0) = P_{lk}(A_{k-l}^c) = 1-p_{lk}.
$$
Therefore,
$$
\tilde{P}_{lk}(\0) = \P(W_{lk}=0) = P_{lk}(\0).
$$
Otherwise, $T=[T_1,\ldots,T_n]\in \T_{k-l-1}^n$ for some $n\geq 1$.
Define $m=|\left\{j:T_j\in A_{k-l-1}\right\}|$. Suppose $m\geq 1$.
Then
\begin{align*}
\tilde{P}_{lk}(T) &= p_{lk} \tilde{Q}_{lk}(T)\\
&=p_{lk}\frac{\P(W_{lk}=n, X_{lk}=m\mid X_{lk}\geq 1)}{{n \choose
m}} \prod_{i:T_i\in A_{k-l-1}} \tilde{Q}_{k-l-1}(T_i)\prod_{j:T_j\in
A_{k-l-1}^c} \tilde{R}_{k-l-1}(T_j)
\end{align*}
Since $m\geq 1$, we have
\begin{align*}
p_{lk}\P(W_{lk}=n, X_{lk}=m\mid X_{lk}\geq 1) &= p_{lk}\frac{\P(W_{lk}=n, X_{lk}=m)}{\P(X_{lk} \geq 1)}\\
&= \P(W_{lk}=n)\P(X_{lk} = m\mid W_{lk}=n)\\
&= \P(W_{lk}=n) {n\choose m}p_{l+1,k}^m(1-p_{l+1,k})^{n-m}.
\end{align*}
Now we use that for $T_i\in A_{k-l-1}$, we have
$\tilde{R}_{l+1,k}(T_i) = 0$ and hence $\tilde{Q}_{l+1,k}(T_i) =
\tilde{P}_{l+1,k}(T_i)/p_{l+1,k}$. Similarly, for $T_j\in A_{k-l-1}^c$,
we have $\tilde{R}_{l+1,k}(T_j) = \tilde{P}_{l+1,k}(T_j)/(1-p_{l+1,k})$.
Plugging this into the expression for $\tilde{P}_{lk}(T)$ and using
the induction hypothesis,
\begin{align*}
\tilde{P}_{lk}(T) = \P(W_{lk}=n)\prod_{i=1}^n {P}_{l+1,k}(T_i) = P_{lk}(T).
\end{align*}
A completely analogous calculation shows that
$\tilde{P}_{lk}(T)=P_{lk}(T)$ if $m=0$. Therefore,
$P_{lk}(T)=\tilde{P}_{lk}(T)$ for all $T\in \T_{k-l}$. \hfill $\Box$

We conclude that if we condition a Galton-Watson tree $\Tb\sim
P_{0k}$ on $A_k$, then $\Tb\sim \tilde{Q}_{0k}$, and if we condition
$\Tb$ on $A_k^c$, then $\Tb\sim \tilde{R}_{0k}$.

\section{The cost of searching a tree}\label{section:treesearch}

We will use the results of the previous section to determine the
expected cost of "searching" a tree. Suppose we draw a random tree
$\Tb\in \T_k$, and we wish to determine if $\Tb$ reaches level
$k$. We start with a searcher in the root at cost 1, and determine
how many children the root has. This will cost us $K$ times the
number of children. Then we move our searcher to one of the
children (chosen at random). This will cost us 1. Now we determine
the number of children at this node of the tree. If at some stage
we get to a node without any children, we backtrack to a node
where there are still undiscovered children. Remember, each time
the searcher moves to a new node, we pay 1, and each time we
determine the number of children of this new node, we pay $K$
times the number of children. We stop when we reach a node at
level $k$. If the tree has no descendants at level $k$, we start
all over with a new tree. However, we do remember the total cost
we made for the unsuccessful tree! What will be the expected cost
for reaching level $k$, given some fixed offspring distribution at level $l$
$W_l\sim \mu_l$?

\subsection{Recursions for the cost}

Fix $k\geq 1$. Define $A_l\subset \T_l$ as the set of all trees of height $l$ that reach level
$l$. Define $C_k$ as the total expected cost of searching a tree until reaching level $k$, including possible restarts. We also define $D_{lk}$ as the expected cost we
make to reach level $k$ when we condition a tree starting at level $l$ to reach level
$k$, and $E_{lk}$ as the expected cost for searching the entire tree
when we condition the tree starting at level $l$ not to reach level $k$. We choose 
\[ D_{kk}=1\ \mbox{and}\ E_{k-1,k}=1.\]
This corresponds to saying that if we reach level $k$, we just pay 1 and we do not look for any children. Furthermore, if we start a tree at level $k-1$ and condition it to not reach level $k$, then we will have $0$ children, and therefore we also only pay 1. We have the
following relations: if $\Tb\in A_k$, then we expect to pay $D_{0k}$,
and if $\Tb\not\in A_k$, we expect to pay $E_{0k}$ and then we start
all over again. The number of unsuccessful trees we have to search
is ${\rm Geo}(p_{0k})$ distributed, so
\[ C_k = (p_{0k}^{-1}-1)E_{0k} + D_{0k}.\]
Furthermore, if we are at level $l$, and we condition on not reaching level $k$, then all children of this tree
will not make it to the last level. The number of
children $\tilde{W}_{lk}$ in this case, according to
\eqref{eq:notinAk}, is distributed according to $\tilde{W}_{lk}\sim
W_{lk} \mid X_{lk}=0$. So
\begin{align*}
E_{lk} &=  1 + \E[K\tilde{W}_{lk} + E_{l+1,k}\tilde{W}_{lk}]\\
&= 1 + (K + E_{l+1,k}) \E[W_{lk}\mid X_{lk}=0].
\end{align*}
Since we know $E_{k-1,k}$, we can now calculate $E_{k-2,k}, \ldots, E_{0k}$. We also have a recursive
relation for $D_{lk}$: our searcher will try the children of a node at level $l$
until she encounters a child that will reach level $k$. Since we are interested in $D_{lk}$, we condition on there being at least one such child. We therefore know that the distribution of $(\tilde{W}_{lk},\tilde{X}_{lk})$ (i.e.,
the total number of children at level $l$ and the number of
children reaching level $k$) in this case is given by
$(\tilde{W}_{lk},\tilde{X}_{lk})\sim (W_{lk},X_{lk})\mid X_{lk}\geq 1$. The
expected number of "dead ends" (children that will not reach level
$k$ of $\T$) tried by the searcher before trying a successful
child, conditioned on $(\tilde{W}_{lk},\tilde{X}_{lk})$, is given by
$(\tilde{W}_{lk} - \tilde{X}_{lk})/(1+\tilde{X}_{lk})$ (there are
$\tilde{W}_{lk} - \tilde{X}_{lk}$ unsuccessful children and each of them
has probability $1/(\tilde{X}_{lk}+1)$ to appear before the first
successful one). Each "dead end" will correspond to an expected
loss of $E_{l+1,k}$. So we get
\begin{align*}
D_{lk} &= 1 + \E\left[K\tilde{W}_{lk} + \frac{\tilde{W}_{lk} - \tilde{X}_{lk}}{1+\tilde{X}_{lk}}\, E_{l+1,k} + D_{l+1,k}\right]\\
&= 1 + D_{l+1,k} + K\E[W_{lk}\mid X_{lk}\geq 1] + E_{l+1,k}\E\left[\frac{W_{lk}
- X_{lk}}{1+ X_{lk}}\mid X_{lk}\geq 1\right].
\end{align*}
We know that $D_{kk}=1$, so this enables us to calculate
$D_{0k}, \ldots, D_{k-1,k}$, and thereby $C_k$.

\subsection{Poisson offspring}

Suppose that all levels have the same Poisson offspring, so $W_{lk} \sim \rm{Pois}(\mu)$. Then also the number of
successful children is Poisson distributed: $X_{lk} \sim
\textrm{Pois}(\mu p_{l+1,k})$. In this case we have $p_{kk}=1$ and for
$l\leq k-1$
\begin{eqnarray*}
1-p_{lk} = \mathbb{P}(X_{lk}=0) = e^{-\mu p_{l+1,k}}.
\end{eqnarray*}
Define $W\sim {\rm Pois}(\mu)$. Elementary calculations show that
\begin{eqnarray*}
\mathbb{E}[W_{lk}\mid X_{lk} = 0] &=&
\frac{\mathbb{E}[W(1-p_{l+1,k})^{W}]}{1-p_{lk}} =
\mu(1-p_{l+1,k}),\\
\mathbb{E}[W_{lk}\mid X_{lk}\geq 1] &=&
\frac{\mu-\mu(1-p_{l+1,k})(1-p_{lk})}{p_{lk}},\\
\mathbb{E}\left[\frac{W_{lk}-X_{lk}}{1+X_{lk}}\mid X_{lk}\geq 1\right] &=&
\frac{1-p_{l+1,k}}{p_{l+1,k}}-\frac{\mu(1-p_{l+1,k})(1-p_{lk})}{p_{lk}}.
\end{eqnarray*}
The resulting recursions for the cost don't seem to be easily
exactly solvable. Nevertheless, they give rise to some interesting
observations. The analysis that follows shows that if $\mu$ is
large, the cost will also be large. The cost also explodes if
$\mu$ is close to zero. Remarkably, in both cases an unsuccessful
tree costs approximately $1$.

\subsubsection{Asymptotics for large $\mu$} If $\mu$ is large,
then $p_{lk}\approx 1$ and $\mu(1-p_{lk})\approx 0$ for all $k$ and $l\leq k$.
Therefore,
\begin{eqnarray*}
E_{lk} &=& 1+(K+E_{l+1,k})\mu(1-p_{l+1,k})\approx 1,\\
D_{lk} &\approx & 1+D_{l+1,k}+K\mu \approx (k-l)(1+K\mu)+D_{kk} \approx (k-l)K\mu,\\
C_k &\approx & D_{0k} \approx kK\mu.
\end{eqnarray*}
This basically means that only the successful tree plays a role in
the cost. Also the costs for ending up in dead ends have no
significant influence on the total search cost. For each level we
just have to pay $K$ times the expected number of children.

\begin{figure}
\includegraphics[width = .95\columnwidth]{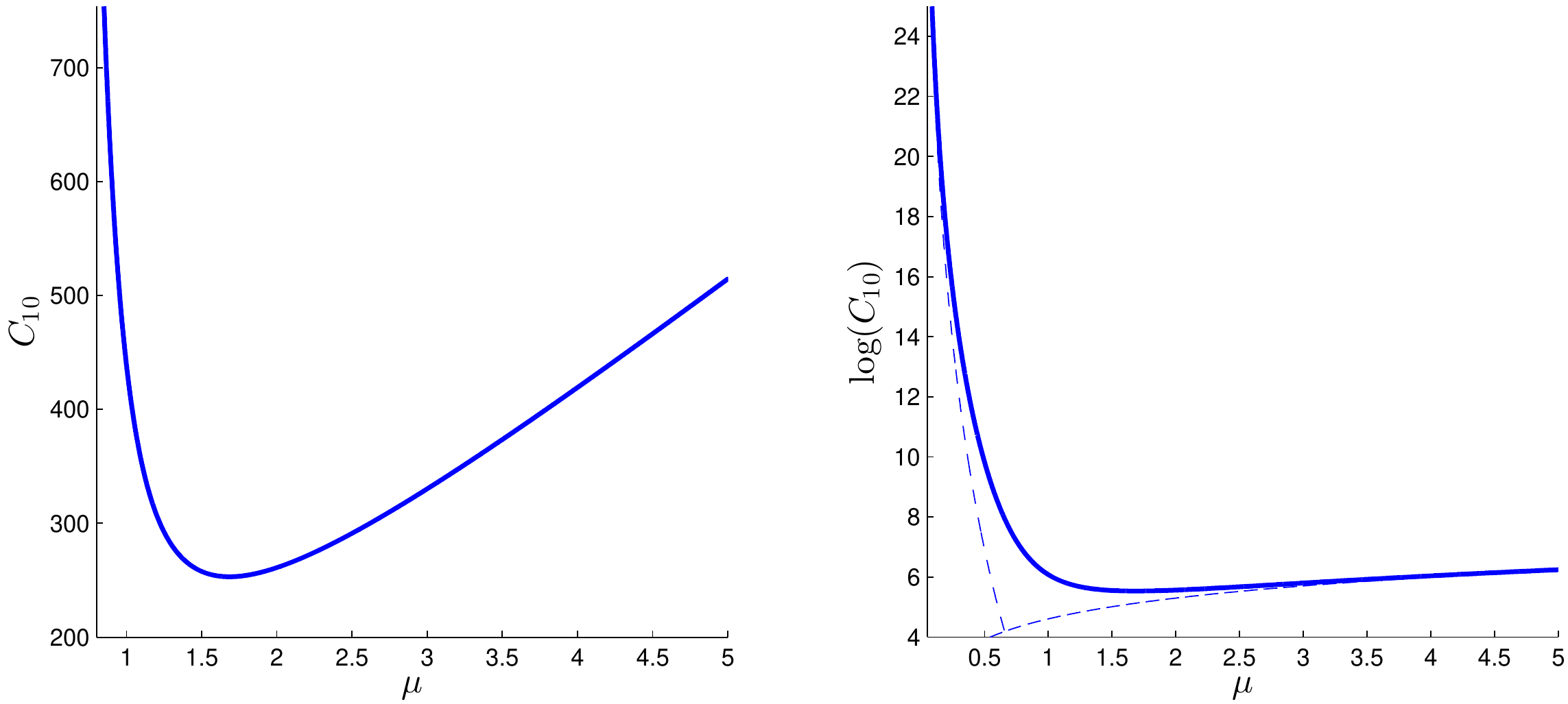}
\caption{Left: plot of the cost as a function of $\mu$ for $K=10$
and $k=10$. Right: plot of the log of the cost together with the
corresponding asymptotic curves.}\label{fig:costplot}
\end{figure}

\begin{figure}
\includegraphics[width = .95\columnwidth]{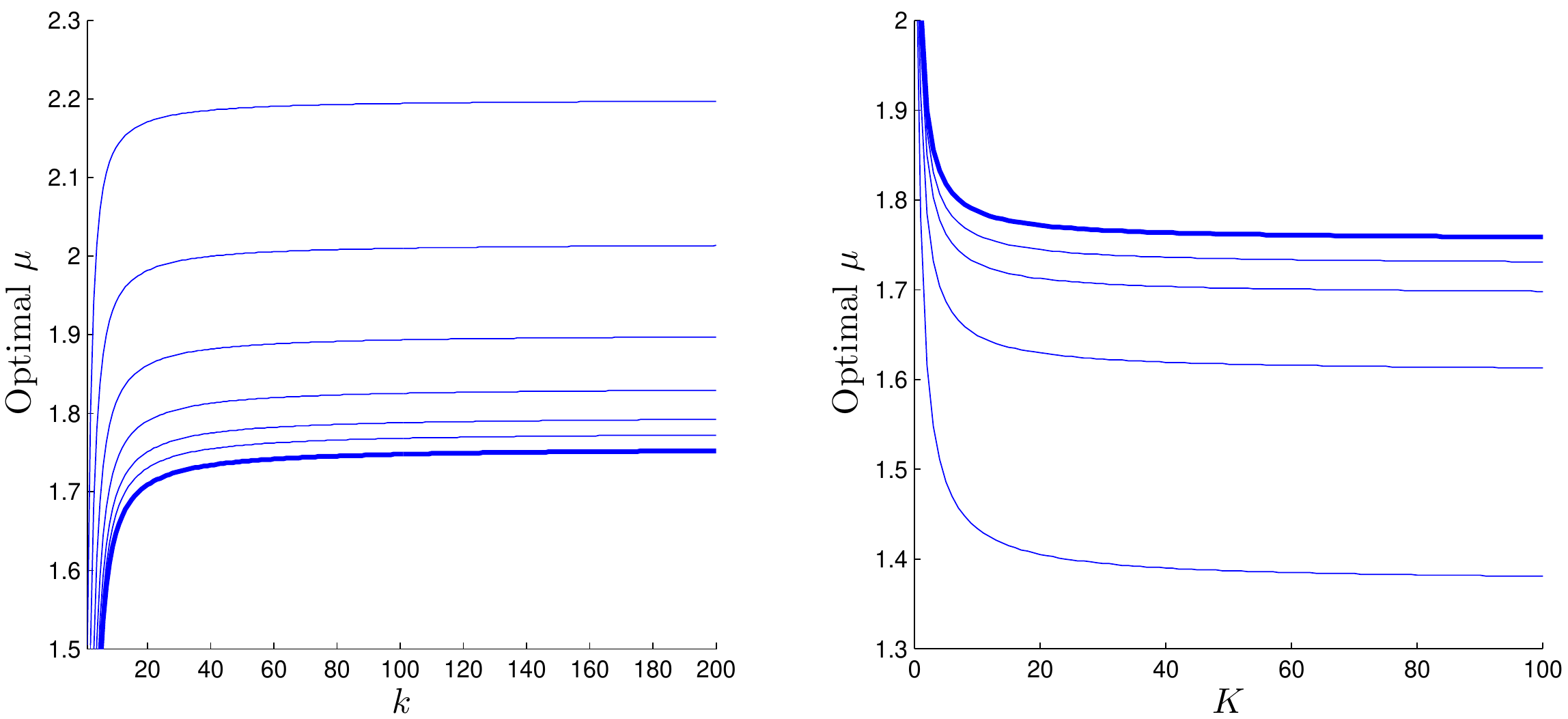}
\caption{Left: the optimal $\mu$ as function of $k$ for
$K=\frac{1}{2},1,2,4,8,16,10^{100}$ (from top to bottom). Right: the optimal $\mu$ as
function of $K$ for $k=4,8,16,32,1000$ (from bottom to top).} \label{fig:optplot}
\end{figure}

\subsubsection{Asymptotics for small $\mu$}

If $\mu$ is close to zero, then
$$
p_{lk} = 1-e^{-p_{l+1,k}\mu} \approx 1-(1-p_{l+1,k}\mu) = p_{l+1,k}\mu,
$$
and since $p_{kk}=1$, we have that $p_{lk}\approx \mu^{k-l}$. It follows
that
\begin{eqnarray*}
E_{lk} &=& 1+(K+E_{l+1,k})\mu(1-p_{l+1,k})\approx  1+K\mu+E_{l+1,k}\mu \approx 1,\\
D_{lk} &\approx& 1+D_{l+1,k}+K\frac{\mu^{k-l}+\mu^{k-l+1}-\mu^{2k-2l}}{\mu^{k-l}}+E_{l+1,k}\frac{\mu^{k-l+1}-\mu^{2k-2l}}{\mu^{k-l}}\\
&\approx & 1+D_{l+1,k}+K+\mu \approx  1+K+D_{l+1,k}\approx
1+k-l+(k-l)K,\\
C_k &\approx & (\frac{1}{\mu^k}-1)E_{0k}+D_{0k} \approx
\frac{1}{\mu^k}+1+k+kK \approx \frac{1}{\mu^k}.
\end{eqnarray*}
Note that here only the unsuccessful trees are important. We
expect $p_{0k}^{-1}$ unsuccessful trees and for each them we pay approximately $1$.

\subsubsection{Numerical results}\label{sec:numerics}

To give an idea of the behaviour of the cost function, we present
some numerical results in this section, without giving any proofs.
Figure \ref{fig:costplot} shows the shape of the cost function for
$K=10$ if the search goal is the $10$th level in the tree. For
comparison we also show the asymptotic curves for small and large
$\mu$ (on a logarithmic scale since the cost rapidly grows as
$\mu$ tends to zero). In this case the cost is minimal when $\mu$
is about $1.68$.

It is interesting to see how the optimal value $\mu_{opt}$ of
$\mu$ changes when we vary $k$ and $K$. The left plot in Figure
\ref{fig:optplot} suggests that for $K$ fixed, $\mu_{opt}$ is
increasing in $k$ and converges to a limit. For small $K$ there is not much to say: if there is no penalty on having many children, $\mu_{opt}$ tends to $\infty$ if $K\rightarrow 0$. Taking $K$ large is more interesting: the right plot in
Figure \ref{fig:optplot} gives evidence that $\mu_{opt}$ also
converges when $k$ is fixed and $K$ goes to $\infty$. The limit
$\lim_{k\rightarrow\infty} \mu_{opt}$ can be interpreted as the
optimal value of $\mu$ in the infinite tree as a function of $K$.
This function is practically the same as the bold curve in the right
panel of Figure \ref{fig:optplot}. A numerical estimate for the
double limit of $\mu_{opt}$ is given by
$$
\lim_{K\rightarrow\infty}\lim_{k\rightarrow\infty}
\mu_{opt}\approx 1.756.
$$

\subsubsection{Search cost in the infinite tree}

In the infinite tree we have two types of nodes: successful ones (i.e. their progeny reaches infinity) and unsuccessful ones. These types do not depend on the level anymore. We define $p$ as the probability that a node is successful. This survival probability satisfies
$$
p = 1-e^{-p\mu}.
$$
We assume that $\mu>1$, to guarantee that $p>0$. The cost for searching an infinite tree is infinite, therefore we consider the cost to `get one step closer to infinity'. More precisely, we look at the expected cost $C_\infty$ to move the searcher from a successful node to a successful child of this node. The number of successful children $X$, given the total number of children $W$ again follows a binomial distribution : $X\mid W\sim {\rm
Bin}(W,p)$. The cost satisfies:
$$
C_\infty = 1+K\cdot\mathbb{E}[W\mid X\geq 1]+E\cdot\mathbb{E}\left[\frac{W-X}{1+X}\mid X\geq 1\right],
$$
where $E$ is the expected cost to search an unsuccessful tree. The second expectation in this expression is the number of unsuccessful trees that have to be searched before picking a successful tree. $E$ satisfies
$$
E = 1+(K+E)\cdot\mathbb{E}[W\mid X=0] = 1+(K+E)\mu(1-p),
$$
and consequently
$$
E = \frac{1+K\mu(1-p)}{1-\mu(1-p)}.
$$
Similar calculations as before lead to
\begin{eqnarray*}
\mathbb{E}[W\mid X\geq 1] &=& \frac{\mu-\mu (1-p)^2}{p},\\
\mathbb{E}\left[\frac{W-X}{1+X}\mid X\geq 1\right] &=& \frac{1-p-\mu(1-p)^2}{p}.
\end{eqnarray*}
Now we evaluate the expected cost:
\begin{eqnarray*}
C_\infty &=& 1+K\cdot\frac{\mu-\mu
(1-p)^2}{p}+\frac{1+K\mu(1-p)}{1-\mu(1-p)}\cdot\frac{1-p-\mu(1-p)^2}{p}\\
&=& \frac{K\mu+1}{p} \approx \frac{\mu}{p}K,
\end{eqnarray*}
where the last approximation holds for $K$ large. Eliminating $\mu$ gives
$$
C_\infty = \frac{-K\log(1-p)+1}{p^2}.
$$
The optimal choice for $\mu$ is therefore
$$
\mu_{opt}(K) = \frac{-\log(1-p_0)}{p_0},
$$
where $p_0$ minimizes $C_\infty$. For $K\rightarrow\infty$, minimizing $\mu/p$ gives (denoting by $W_{-1}$ the lower branch of the real-valued Lambert $W$ function)
$$
\mu_{opt} = \frac{\frac{1}{2}+W_{-1}\left(\frac{-1}{2\sqrt{e}}\right)}{e^{\frac{1}{2}+W_{-1}\left(\frac{-1}{2\sqrt{e}}\right)}-1} \approx 1.756,
$$
which is equal to the numerical estimate obtained in Section \ref{sec:numerics}.

\section{Generalization}\label{section:generalization}

So far we have focused on trees conditioned on survival to a fixed
level, but the ideas can be easily extended to a more general
setting. As a first simple example, we can define $A_l$ to be the set of all
trees $[T_1,\ldots,T_n]\in\T_{l-1}^n$ such that at least
\emph{two} of the children are in $A_{l-1}$. The measures
$\tilde{Q}_{lk}$ and $\tilde{R}_{lk}$ will again be defined as in
(\ref{eq:inAk}) and (\ref{eq:notinAk}) for trees in $A_{k-l}$ and in
$A_{k-l}^c$ respectively. Only the conditioning on $X_{lk}$ in both
equations has to be adapted. The resulting measure $\tilde{P}_{0k}$ can be used to construct a tree conditioned on having a full binary subtree reaching level $k$.

In the examples discussed so far, the set $\T_l$ is partitioned into two subsets: $A_l$ and $A_l^c$. To state it differently: we distinguished between two types of trees, trees that are in some sense successful and those that are not. Allowing for partitions into more than two subsets widens the spectrum of events we can condition on and, as we will show, the machinery still works. We will discuss two examples where it is natural to consider three types of trees. In the first example, we condition a tree of height $k$ on the property that each node (except those at the last two levels) has at least two grandchildren. In the second example, we condition the tree on having height exactly $k$, thus producing an alternative for the construction of Geiger and Kersting \cite{Geiger}. We will now set up our general framework and show how these examples fit into it.       

We start by choosing $k_0\in\N$ and partitioning $\T_{k_0}$ into $m$ subsets $A_{k_0}^{(1)},\ldots,A_{k_0}^{(m)}$. This partition will be the starting point to recursively define partitions of $\T_l,l>k_0$ into sets $A_{l}^{(i)}, i=1,\ldots,m$. Fix $k>k_0$: this will be the (maximum) size of the trees we are considering. We adopt the notation $\tilde T \prec T$ to state that $\tilde T$ is a child of $T$. For trees in $\T_l,l>k_0$, we recursively define a function $N_l:\T_l\rightarrow\N^m$ that counts how many children of each type an element $T$ of $\T_l$ has: 
$$
N_l(T) = \left(N_{l,1}(T),\ldots,N_{l,m}(T)\right),\quad\textrm{where}\quad N_{l,i}(T) = \#\left\{\tilde T\prec T:\tilde T\in A_{l-1}^{(i)}\right\}.  
$$ 
Next, for each $l>k_0$, we partition $\N^m$ into subsets $B_{l,1},\ldots,B_{l,m}$. This partition is the key for the recursive definition of $A_{l}^{(i)}$. The set $A_{l}^{(i)}$ will contain exactly those trees for which the counting vector $N_l(T)$ is in $B_{l,i}$:
$$
A_{l}^{(i)} = \left\{T\in\T_l:N_l(T)\in B_{l,i}\right\}.
$$    

\subsection{Examples}
Let's see how the examples mentioned in this section fit into this framework. We start with a tree conditioned to have a full binary tree up to level $k$. Choose $k_0=1$ and partition $\T_1$ into the two sets 
\[A_1^{(1)}=\{ T\in \T_1\mid \#\{\tilde{T}\prec T\}\geq 2\} \mbox{  and  } A_1^{(2)}=\{ T\in \T_1\mid \#\{\tilde{T}\prec T\}\leq 1\}.\]
Then define for each $l>1$ 
\[ B_{l,1} = \{ (n_1,n_2)\in\N^2\mid n_1\geq 2\}\mbox{  and  } B_{l,2}=\N^2\setminus B_{l,1}.\]
For $T\in\T_k$, one easily checks that $T\in A^{(1)}_k$ precisely when $T$ has a full binary subtree reaching level $k$.

As a second example, consider trees in $T_k$ where each node (except at the last two levels) has at least two grandchildren. We start by choosing $k_0=2$ and partitioning $\T_2$ into three sets:
\begin{align*}
 A_2^{(1)} & =\{ T\in \T_2\mid T\mbox{ has one child, at least two grandchildren}\},\\
A_2^{(2)} & =\{T\in \T_2\mid T\mbox{ has at least two children and at least two grandchildren}\},\\
A_2^{(3)} & = \{T\in \T_2\mid T\mbox{ has at most one grandchild}\}.
\end{align*}
Define for each $l>2$
\[ B_{l,1} = \{(0,1,0)\}, B_{l,2}=\{ n\in \N^3\mid n_1+n_2\geq 2\}, B_{l,3}=\N^3\setminus(B_{l,1}\cup B_{l,2}).\]
This construction guarantees that if a tree $T\in \T_k$ is of type 1, then the root has precisely one child, which is of type 2. Furthermore, $T$ is of type 2 if the root has at least 2 children of either type 1 or type 2. Otherwise, $T$ is of type three. The set $A_k^{(1)}\cup A_k^{(2)}$ contains exacty all trees in $\T_k$ such that each node up to level $k-2$ has at least two grandchildren. We need three types in this construction, because if we only check that the children of a node have at least two grandchildren, then this node itself may have only one grandchild (a node is allowed to have only one child).

As a last example, consider a tree in $\T_k$ that is conditioned to reach level $k-1$, but not level $k$. We start by choosing $k_0=2$, and partitioning $\T_2$ into three sets, namely correct trees, short trees and long trees:
\begin{align*}
 A_2^{(1)} & =\{ T\in \T_2\mid T\mbox{ reaches level 1, but not level 2}\},\\
A_2^{(2)} & =\{T\in \T_2\mid T\mbox{ does not reach level 1}\} = \{\0\},\\
A_2^{(3)} & = \{T\in \T_2\mid T\mbox{ reaches level 2}\}.
\end{align*}
Define for each $l>2$
\begin{align*}
 B_{l,1} & = \{ n\in\N^3\mid n_1\geq 1, n_3=0\},\\
B_{l,2} & = \{ n\in\N^3\mid n_1=0, n_3=0\},\\
B_{l,3} & = \{ n\in\N^3\mid n_3\geq 1\},\\
\end{align*}
This construction guarantees that if a tree $T\in \T_k$ is of type 1, then it has at least one child that reaches level $k-1$, and no children that reach level $k$. If $T$ is of type 2, all its children do not reach level $k-1$, and if $T$ is of type three, then at least one child reaches level $k$. Conditioning on being in $A_k^{(1)}$ therefore gives the desired result.

\subsection{Conditional measures}

Define $p_{lk}^{(i)}$ for $0\leq l\leq k-k_0$ by
$$
p_{lk}^{(i)} = P_{lk}(A_{k-l}^{(i)}).
$$
We can calculate this probability in a recursive way. Denote, for $q\in [0,1]^m$ with $\sum q_i =1$, by ${\rm Multi}(n,q)$ the multinomial distribution where we distribute $n$ elements over $m$ types, according to the probabilities $q_i$. We also choose independent random variables $W_{lk}\sim \mu_l$ according to the offspring distribution at level $l$. Then, for $l<k-k_0$
\[ p_{lk}^{(i)} = \P\left( {\rm Multi}\left(W_{lk}, (p_{l+1,k}^{(1)}, \ldots, p_{l+1,k}^{(m)})\right) \in B_{l,i}\right).\]
As before, we want to define the conditional measure $\tilde{Q}_{lk}^{(i)}$ on $A_{k-l}^{(i)}$. To do this, define the independent random vectors $X_{lk} =  (X_{lk}^{(1)}, \ldots, X_{lk}^{(m)})$ such that
\[ (X_{lk}^{(1)}, \ldots, X_{lk}^{(m)})\mid W_{lk}\sim {\rm Multi}\left(W_{lk}, (p_{l+1,k}^{(1)}, \ldots, p_{l+1,k}^{(m)})\right).\]
Of course, this represents the number of children in each of the types.
For each $N\in \N^m$ we define the combinatorial constant
\[ D(N) = \frac{(\sum_{i=1}^m N_i)!}{N_1!\cdots N_m!}.\]
For $l=k-k_0$, we define for each $T\in A^{(i)}_{k_0}$
\[ \tilde{Q}_{k-k_0,k}^{(i)}(T) = \frac{P_{k-k_0,k}(T)}{P_{k-k_0,k}(A^{(i)}_{k_0})}.\]
Next, we inductively define the measures $\tilde{Q}_{lk}^{(i)}$ on $A_{k-l}^{(i)}$ for each $0\leq l\leq k-k_0-1$ such that for each $T\in A^{(i)}_{k-l}$
\[ \tilde{Q}^{(i)}_{lk}(T) = \frac{\P(X_{lk}= N_{k-l}(T)\mid X_{lk}\in B^{(i)}_{k-l})}{D(N_{k-l}(T))}\prod_{j=1}^m\prod_{\tilde{T}\prec T : \tilde{T}\in A_{k-l-1}^{(j)}} \tilde{Q}^{(j)}_{l+1,k}(\tilde{T}).\]
As usual, empty products are taken to be 1. Note that this definition is valid for all $T\in\T_{k-l}$: we simply get $\tilde{Q}_{lk}^{(i)}(T)=0$ whenever $T\not\in A_{k-l}^{(i)}$. As before, we can now define the alternative measure $\tilde{P}_{lk}$ on $\T_{k-l}$:
\[ \tilde{P}_{lk}(T)= \sum_{i=1}^m p_{lk}^{(i)}\tilde{Q}_{lk}^{(i)}(T).\]

\begin{stelling} For all $0\leq l\leq k-k_0$ and $T\in \T_{k-l}$,
$$
P_{lk}(T)=\tilde{P}_{lk}(T).
$$
\end{stelling}
{\bf Proof:} The theorem is true by construction for $l=k-k_0$. Now suppose that we have already shown that $P_{l+1,k}=\tilde{P}_{l+1,k}$. Choose $T\in \T_{k-l}$ and suppose $T\in A_{k-l}^{(i)}$.  
\begin{align*}
\tilde{P}_{lk}(T) & = p_{lk}^{(i)}\tilde{Q}_{lk}^{(i)}(T)\\
&= \frac{p_{lk}^{(i)}\P(X_{lk}= N_{k-l}(T)\mid X_{lk}\in B^{(i)}_{k-l})}{D(N_{k-l}(T))}\prod_{j=1}^m\prod_{\tilde{T}\prec T : \tilde{T}\in A_{k-l-1}^{(j)}} \tilde{Q}^{(j)}_{l+1,k}(\tilde{T})\\
& = \frac{p_{lk}^{(i)}\P(W_{lk}=\sum_{j=1}^m N_{k-l,j}(T))}{\P(X_{lk}\in B_{k-l}^{(i)})}\left(\prod_{j=1}^m \left(p_{l+1,k}^{(j)}\right)^{N_{k-l,j}(T)}\right)\prod_{j=1}^m\prod_{\tilde{T}\prec T : \tilde{T}\in A_{k-l-1}^{(j)}} \tilde{Q}^{(j)}_{l+1,k}(\tilde{T}).
\end{align*}
Note that $\#\{\tilde{T}\prec T\mid \tilde{T}\in A^{(j)}_{k-l-1}\}=N_{k-l,j}(T)$, so we get
\begin{align*}
\tilde{P}_{lk}(T) & = \frac{p_{lk}^{(i)}\P(W_{lk}=\#\{\tilde T\prec T\})}{\P(X_{lk}\in B_{k-l}^{(i)})}\prod_{j=1}^m\prod_{\tilde{T}\prec T : \tilde{T}\in A_{k-l-1}^{(j)}} p_{l+1,k}^{(j)}\tilde{Q}^{(j)}_{l+1,k}(\tilde{T})\\
& = \frac{p_{lk}^{(i)}\P(W_{lk}=\#\{\tilde T\prec T\})}{\P(X_{lk}\in B_{k-l}^{(i)})}\prod_{\tilde T\prec T} \tilde{P}_{l+1,k}(\tilde{T})\\
& = \frac{p_{lk}^{(i)}\P(W_{lk}=\#\{\tilde T\prec T\})}{\P(X_{lk}\in B_{k-l}^{(i)})}\prod_{\tilde T\prec T} P_{l+1,k}(\tilde{T})\\
& = \frac{p_{lk}^{(i)}}{\P(X_{lk}\in B_{k-l}^{(i)})}\,P_{lk}(T).
\end{align*}
As a final step we use that
\begin{align*}
\P(X_{lk}\in B_{k-l}^{(i)}) & = \sum_{n=0}^\infty \P(W_{lk}=n)\P\left(  {\rm Multi}\left(n, (p_{l+1,k}^{(1)}, \ldots, p_{l+1,k}^{(m)})\right)\in B_{k-l}^{(i)}\right), 
\end{align*}
and this is exactly the probability that a tree under $P_{lk}$ is an element of $A_{k-l}^{(i)}$, since $P_{l+1,k}(A_{k-l-1}^{(j)}) = p_{l+1,k}^{(j)}$. This proves that $\P(X_{lk}\in B_{k-l}^{(i)}) = p_{lk}^{(i)}$, which implies that $\tilde{P}_{lk}(T)=P_{lk}(T)$.
\hfill $\Box$

As before, conditioning $\T\sim P_{lk}$ on $\T\in A_{k-l}^{(i)}$ simply means that the tree follows $\tilde{Q}_{lk}^{(i)}$.

\end{document}